\documentclass{ifacconf}

\usepackage{graphicx}      
\usepackage{natbib}        
\usepackage{amsmath}

\usepackage{amsmath} 
\usepackage{amssymb}
\usepackage{amsfonts,dsfont, bm, eufrak} 
\usepackage{color}

\usepackage{tikz} 
\usepackage{circuitikz}
\usepackage[customcolors]{hf-tikz}
\usetikzlibrary{arrows,shapes,backgrounds,decorations.text}

\usetikzlibrary{arrows.meta}

\usetikzlibrary{backgrounds}
\pgfdeclarelayer{myback}
\pgfsetlayers{background,myback,main}

\usepackage{datetime}
\usepackage{enumitem}

\usepackage{graphicx}
\graphicspath{{./img/},{./fig/},{./FiguresFB/},{./fig_time_var_sim_2/} }
\newlength{\figtriml}
\newlength{\figtrimb}
\newlength{\figtrimr}
\newlength{\figtrimt}
\newlength{\figwidth}
\newlength{\figtrimlFB}
\newlength{\figtrimbFB}
\newlength{\figtrimrFB}
\newlength{\figtrimtFB}
\newlength{\figtrimlFBmat}
\newlength{\figtrimlFBX}
\newlength{\figtrimbFBX}
\newlength{\figtrimrFBX}
\newlength{\figtrimtFBX}
\newlength{\figtrimlW}
\newlength{\figtrimbW}
\newlength{\figtrimrW}
\newlength{\figtrimtW}
\newcommand{\field}{\mathfrak{f}}

\renewcommand{\l}{\left}
\renewcommand{\r}{\right}

\DeclareMathOperator*{\diag}{Diag}

\newcommand{\1}{\mathbf{1}}
\newcommand{\0}{\mathbf{0}}

\newcommand{\Rp}{\mathbb{R}_+}

\newcommand{\GG}{\mathcal{G}}

\newcommand{\xx}{\mathbf{x}}

\newcommand{\Wil}{W^{i\to \ell}}
\newcommand{\Hil}{H^{i\to \ell}}
\newcommand{\Wij}{W^{i\to j}}
\newcommand{\Hij}{H^{i\to j}}
\newcommand{\Wki}{W^{k\to i}}
\newcommand{\Hki}{H^{k\to i}}

\newcommand{\Wfj}{W^{\field\to j}}
\newcommand{\Hfj}{H^{\field\to j}}

\newcommand{\bWij}{\bar{W}^{i\to j}}
\newcommand{\bHij}{\bar{H}^{i\to j}}

\newtheorem{defi*}{Definition}

\begin{document}
\begin{frontmatter}

\title{Effects of Network Communities and Topology Changes in Message-Passing Computation of Harmonic Influence in Social Networks}



\author[First]{W. S. Rossi} 
\author[Second]{P. Frasca} 

\address[First]{University of Twente, 7500 AE Enschede, The Netherlands (e-mail: w.s.rossi@utwente.nl).}
\address[Second]{Univ.\ Grenoble Alpes, CNRS, Inria, Grenoble INP, GIPSA-lab, 38000 Grenoble, France (e-mail: paolo.frasca@gipsa-lab.fr).}


\begin{abstract}                
The harmonic influence is a measure of the importance of nodes in social networks, which can be approximately computed by a distributed message-passing algorithm. 
In this extended abstract we look at two open questions about this algorithm. How does it perform on real social networks, which have complex topologies structured in communities? How does it perform when the network topology changes while the algorithm is running? 
We answer these two questions by numerical experiments on a Facebook ego network and on synthetic networks, respectively.
We find out that communities can introduce artefacts in the final approximation and cause the algorithm to overestimate the importance of ``local leaders'' within communities.
We also observe that the algorithm is able to adapt smoothly to changes in the topology.
%
\end{abstract}

\begin{keyword}
Distributed algorithm, Message-passing, Opinion dynamics, Social networks 
\end{keyword}

\end{frontmatter}

\section{Harmonic Influence and Message Passing}

In the study of social networks and dynamical processes therein, it is important to identify the most influential leaders. Several definitions have been used to evaluate nodes as potential leaders, e.g.
~\cite{linfarjovTAC14leaderselection,7277027,Mieghem:2017:spreader}. 
The \textit{harmonic influence} is a definition motivated by a linear opinion dynamics model with stubborn agents. It was introduced in~\cite{Vassio:2014:journal} and implicitly used in~\cite{como:opinion,Yildiz:2013}. 
We recall its equivalent definition given by~\cite{RF:2016:TNSE}. 
Consider a simple weighted graph\footnote{Vectors are denoted with boldface letters and matrices with capital letters. 
The all-zero and all-one vectors are denoted by $\0$ and $\1$, respectively. A graph is said to be connected if for any pair of nodes $i$, $j$ there exists a sequence of adjacent edges that joins them.} $\GG = (I,E,C)$ with node set $I= \{\field,1,2,\ldots ,n\}$ where $\field$ is a special node called \textit{field}. 
The edge set $E$ contains unordered pairs of nodes and the non-negative weight matrix $C \in \Rp^{I\times I}$ is such that $C_{ij}$ and $C_{ji}$ are both non-zero if and only if $\{i,j\} \in E$. 
%
%
%
We also introduce the diagonal matrix $D = \diag(C \1)$  and the {\em Laplacian} matrix $L = D -C\,.$
We assume $C$ to be symmetric and the graph $\GG$ to be connected.
%
Given a node $\ell \neq \field$ where $\ell$ stands for \textit{leader}, let $R^\ell := I\setminus \{\field,\ell\}$ be the set of remaining nodes and consider the discrete Dirichlet problem 
\begin{align} \label{eq:lap-sys} 
	\left\{ \begin{array}{l} 
	\l( L \, \xx \r)_{R^\ell}= \0 \\ 
	x_\ell = 1 \\
	x_\field = 0 \,. \end{array}\right. 
\end{align}
The \textit{harmonic influence} of $\ell$ is the sum of entries of the vector $\xx$ solution of~\eqref{eq:lap-sys}, that is,
\begin{align}\label{eq:def-H-ell}
	H(\ell) := \1^\top \xx\,.
\end{align}
%
Then, the computation of the harmonic influence of the $n$ nodes in  $I \setminus \{ \field\}$ requires the solution of  $n$ linear systems. A naive approach would then require global knowledge of the graph and would not exploit apparent redundancies in the computations. 
To overcome these issues, \cite{Vassio:2014:journal} proposed the following distributed \textit{Message Passing Algorithm} (MPA) that computes the influences of all nodes at the same time. 
%

Let $t \in\{0,1,\ldots\}$ be an  iteration counter and let the set $N_i = \{j \in I: \{i,j\}\in E \}$ contain the neighbors of $i$ in $\GG$.
At each step, every node $i$ sends to all its neighbors $j$ two messages:
$$\Wij(t) \in [0,1]\,, \qquad \Hij(t) \in [0,+\infty) \,.$$
The field node $\field$ sends null messages:
$$\Wfj(t)=0\,, \quad \Hfj(t)=0\,,\quad \forall j\in N_\field\,,\quad \forall  t\geq0\,,$$
whereas any other node $i\neq \field$ sends the initial messages: 
\begin{align}\label{eq:MPA-initial}	\Wij(0) = 1\,,\quad  \Hij(0) = 1\,,\quad \forall j\in N_i\,\end{align}
and then synchronously updates the messages sent to his neighbor $j$ following the rules:
\begin{align}
\label{eq:MPA-update-1-vassio}
&\Wij(t+1) = \left( 1 + { \textstyle{\sum_{k \in N_i^j }} \frac{C_{ik}}{C_{ij}}} \l(1- \Wki(t) \r)  \right)^{-1} \\  
&\Hij(t+1) = 1 + { \textstyle{\sum_{k \in N_i^j }} }\Wki(t)\, \Hki(t)\,, \label{eq:MPA-update-2} 
\end{align} 
where $N_i^j:= N_i \setminus \{j\}$ is the set of neighbors of $i$ except the one to which the message is sent. At any time, any node $\ell$  in $I\setminus \{\field\} $ can compute an approximation of 
$H(\ell)$ by
\begin{align*}
	H^\ell(t)  =  1 + {\textstyle{ \sum_{i \in N_\ell }} }\Wil(t)\, \Hil(t)\,. \nonumber
\end{align*}

The MPA is exact on trees, where it converges in a number of steps equal to the diameter of the graph. 
On general graphs, the algorithm converges asymptotically as proved in~\cite{RF:2016:ecc,RF:2016:TNSE}. Based on extended simulations on random graphs, the typical convergence time of the algorithm is conjectured to be $O(m/n)$, where $m$ is the number of edges. A mean-field argument by ~\cite{RF:2017:ifac} corroborates this conjecture for homogeneous networks.
In general, the limit values overestimate the exact values of the harmonic influence (that is, $H^\ell(\infty) \geq H(\ell)$). 
However, on random graphs the ranking between the nodes that is provided by the algorithm 
is in very good agreement with the exact ranking. 
%
%
%
%
%
%
%

In Section~\ref{sec:COMM}, we look at the correctness and convergence time for real social networks, which have a pronounced community structure. While convergence is guaranteed by the theoretical results, it not obvious whether that the convergence time follows the conjecture and whether the ranking remains meaningful. 
In Section~\ref{sec:T-VAR}, we extend the algorithm to networks that change topology while the MPA computation unfolds. For this case, convergence is not guaranteed by the available theory.

\section{The Effects of Community Structure}\label{sec:COMM}

Real social networks have complex topologies of interconnections that are often organized in communities.
On one side, the degree distribution of real networks (i.e. the distribution of the number of interconnection of each node) is typically  broad, with relatively few nodes of high degree and many of low degree. On the other side, nodes can be grouped in communities, such that most of the edges are concentrated within each community \cite{Fortunato:2010:comm}. 
%
%
In this section, we test the MPA on a real social network with community structure, extracted from the dataset\footnote{http://snap.stanford.edu/data/egonets-Facebook.html} collected by~\cite{LM:2012:NIPS}. 

\begin{figure}	
\centering 	
\fbox{%
	\includegraphics[trim={\figtrimlFBmat} {\figtrimbFB} {\figtrimrFB} {\figtrimtFB},
	clip, width={\figwidth}, keepaspectratio=true]{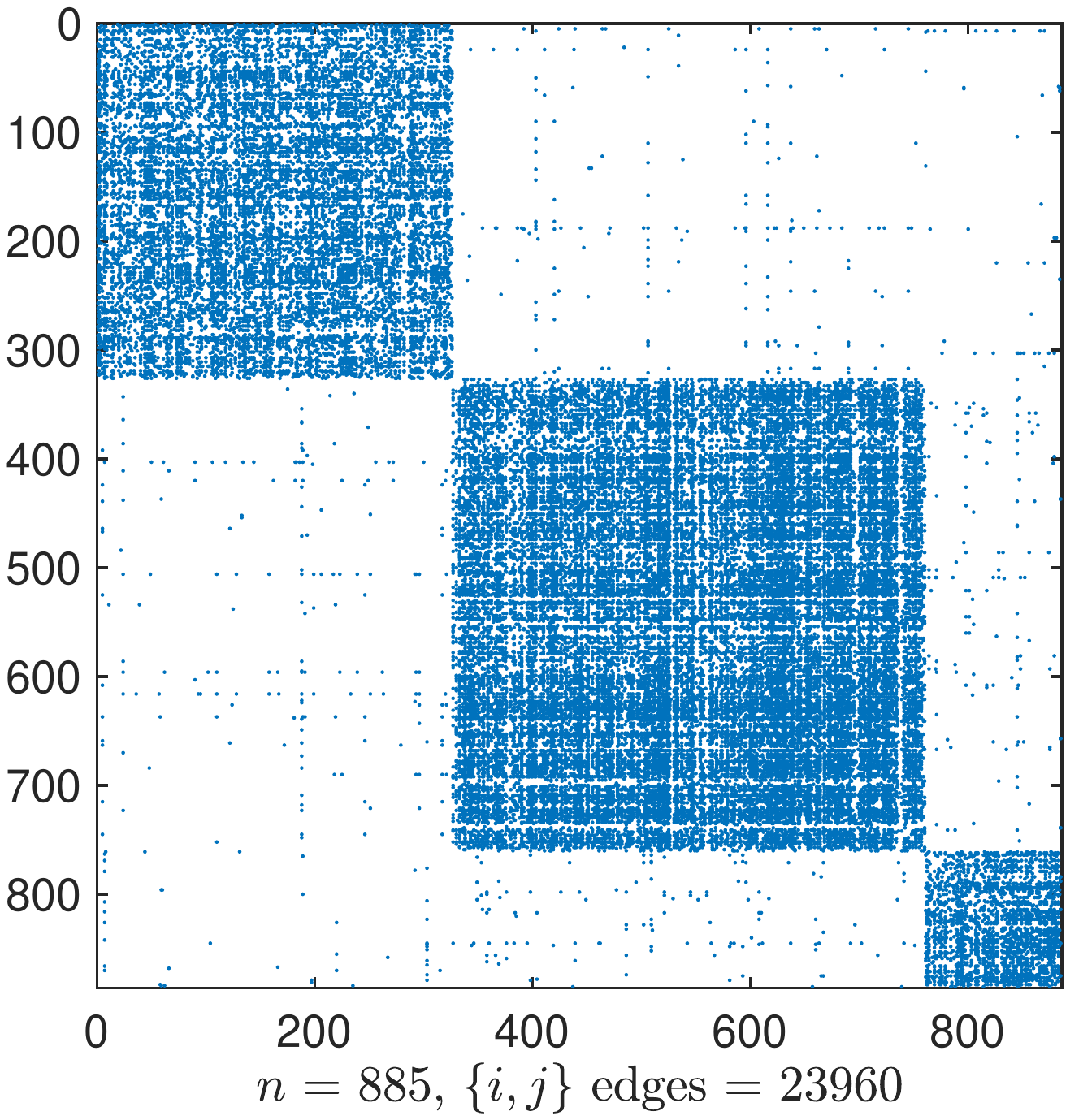}%
	}
	\caption{\label{fig:comm_mm_1_no_seed_spy} The adjacency matrix of the selected ego network. The rows and columns are reordered according to the community structure, showing three prominent communities with 326, 434 and 125 nodes, respectively. A closer inspection reveals that within each communities coexist nodes with high and low degree.	}
\end{figure}

The dataset contains a collection of \textit{ego networks} from the Facebook social graph. Let $F = (U,A)$ be the full unweighed Facebook social graph at the data collection time, where $U$ is the users set and  $A$ the edge set representing the  acquaintance relations. Given a user $u \in U$, the set $N_u$ is the set of Facebook friends of $u$. The ego network of $u$ is the subgraph of $F$ induced by $N_u$, i.e 
$$\left(N_u, \big\{ \{v,w\}  \in A : v,w\in N_u\big\} \right)\,,$$
and does not contain $u$. 
From the dataset we extracted one ego network with 885 nodes and 23960 edges. 
We identified the communities of the ego network with the ``Louvain algorithm'' by \cite{Blondel:2008:algo-modularity}; in particular, we used the implementation\footnote{https://sites.google.com/site/bctnet} by \cite{Rubinov:2010:BCT-toolbox}  with default parameters. We found three prominent communities with 326, 434 and 125 nodes, see Fig.~\ref{fig:comm_mm_1_no_seed_spy}. 

From the ego network we build two weighted graphs for the simulations. 
%
The graph $\GG_1 = (I_1,E_1,C_1)$ contains all the ego network augmented by the field node. The node set is $I_1 = \{\field,1,\ldots,885\}$ and the edge set $E_1$ contains all the 23960 edges of the ego network. 
The edge set also include every edge of the form $\{\field, i\}$ with $i\in \{1,\ldots,885\}$: we can interpret the field node $\field$ as the original user $u$ and these edges as those between $u$ and his friends in $N_u$.
Finally, the entries of the matrix $C_1 \in \Rp^{I_1\times I_1}$ are
\begin{equation*}
\left\{ \begin{array}{ll} 
	(C_1)_{i\field}= (C_1)_{\field i} = 0.040  & \text{ for every } i \in \{1,\ldots,885\} \\ 
	(C_1)_{ij}= 1 & \text{ if } i,j\neq \field \text{ and } \{i,j\}\in E_1 \\
	(C_1)_{ij}= 0 & \text{ if } \{i,j\}\notin E_1 \end{array}\right. 
\end{equation*}
The graph $\GG_2 = (I_2,E_2,C_2)$ is the subgraph of $\GG_1$ induced by $I_2 = \{\field,327,\ldots,760\}$. 
It is restricted to the second community of the ego network, augmented with the field node and his edges. 
The edge set $E_2$ contains 16253 edges between the 434 non-field nodes; the matrix $C_2$ follows accordingly.

\begin{figure}	
\centering 	
\fbox{%
	\includegraphics[trim={\figtrimlFB} {\figtrimbFBX} {\figtrimrFB} {\figtrimtFBX},
	clip, width={\figwidth}, keepaspectratio=true]{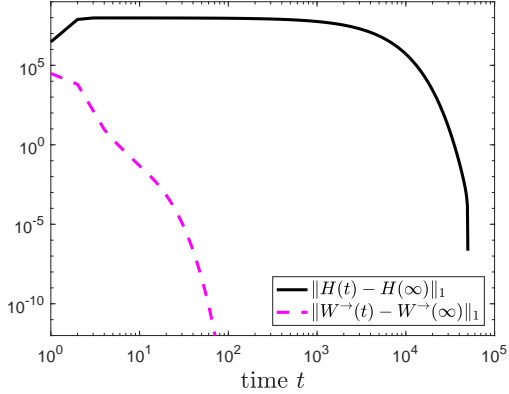}%
	}
	\caption{\label{fig:comm_no_seed_5_dyn} The convergence of the MPA on the network $\GG_2$. The solid black line is the distance to convergence of the estimates of the harmonic influence obtained by the MPA. The dashed magenta line is the distance to convergence of the messages $\Wij(t)$.}
\end{figure}

\begin{figure}	
\centering 	
\fbox{%
	\includegraphics[trim={\figtrimlFB} {\figtrimbFB} {\figtrimrFB} {\figtrimtFB},
	clip, width={\figwidth}, keepaspectratio=true]{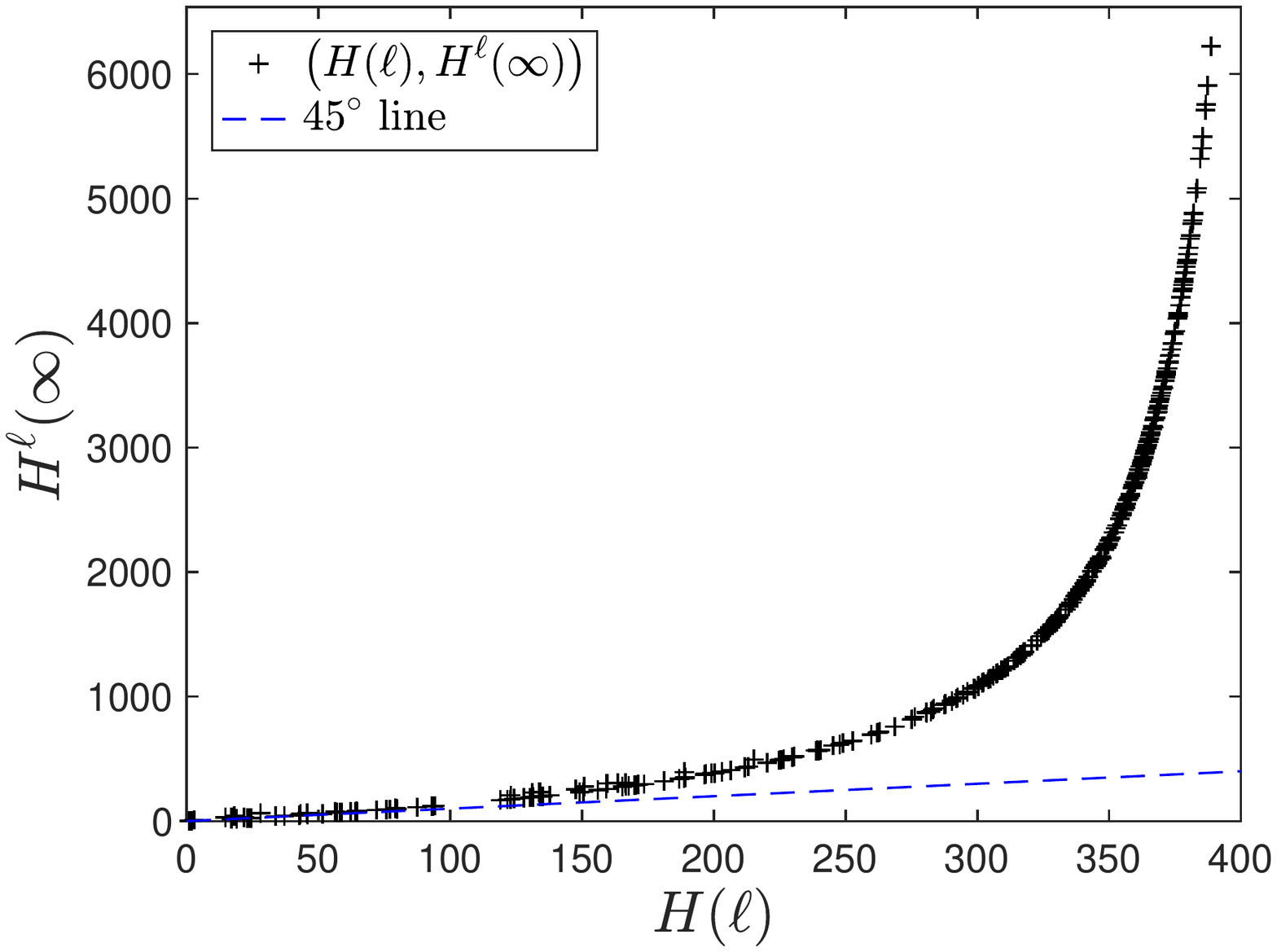}%
	}
	\caption{\label{fig:comm_no_seed_5_HH} The asymptotic values $H^\ell(\infty)$ of the harmonic influence computed by the MPA against the corresponding exact values $H(\ell)$ computed by the definition, for the graph $\GG_2$. All crosses are above the $45^\circ$ line. }
\end{figure}

We first discuss the simulation on the graph $\GG_2$, with a single community but with nodes of very different degree. 
Fig.~\ref{fig:comm_no_seed_5_dyn} represents the convergence of the MPA: the $\Wij(t)$ messages take about 20 steps to converge while the estimates $H^\ell(t)$  of the harmonic influence are very slow, taking almost 49000 iterations. 
Fig.~\ref{fig:comm_no_seed_5_HH} compares the estimates $H^\ell(\infty)$ with the exact values $H(\ell)$ computed with the definition. The MPA algorithm largely overestimates the harmonic influence, but the ranking remains well preserved. 
 

\begin{figure}	
\centering 	
\fbox{%
	\includegraphics[trim={\figtrimlFB} {\figtrimbFBX} {\figtrimrFB} {\figtrimtFBX},
	clip, width={\figwidth}, keepaspectratio=true]{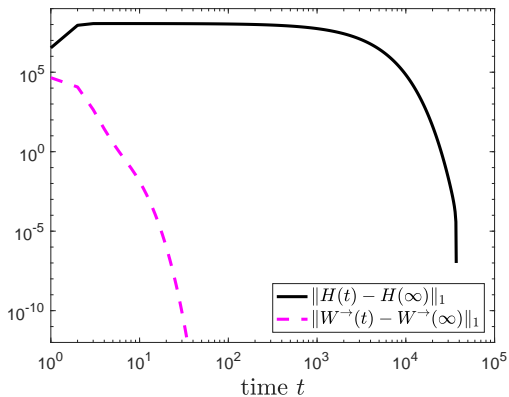}%
	}
	\caption{\label{fig:comm_mm_1_no_seed_dyn} The convergence of the MPA on the network $\GG_1$. The solid black line is the distance to convergence of the estimates of the harmonic influence obtained by the MPA. The dashed magenta line is the distance to convergence of the messages $\Wij(t)$. }
\end{figure}

\begin{figure}	
\centering 	
\fbox{%
	\includegraphics[trim={\figtrimlFB} {\figtrimbFB} {\figtrimrFB} {\figtrimtFB},
	clip, width={\figwidth}, keepaspectratio=true]{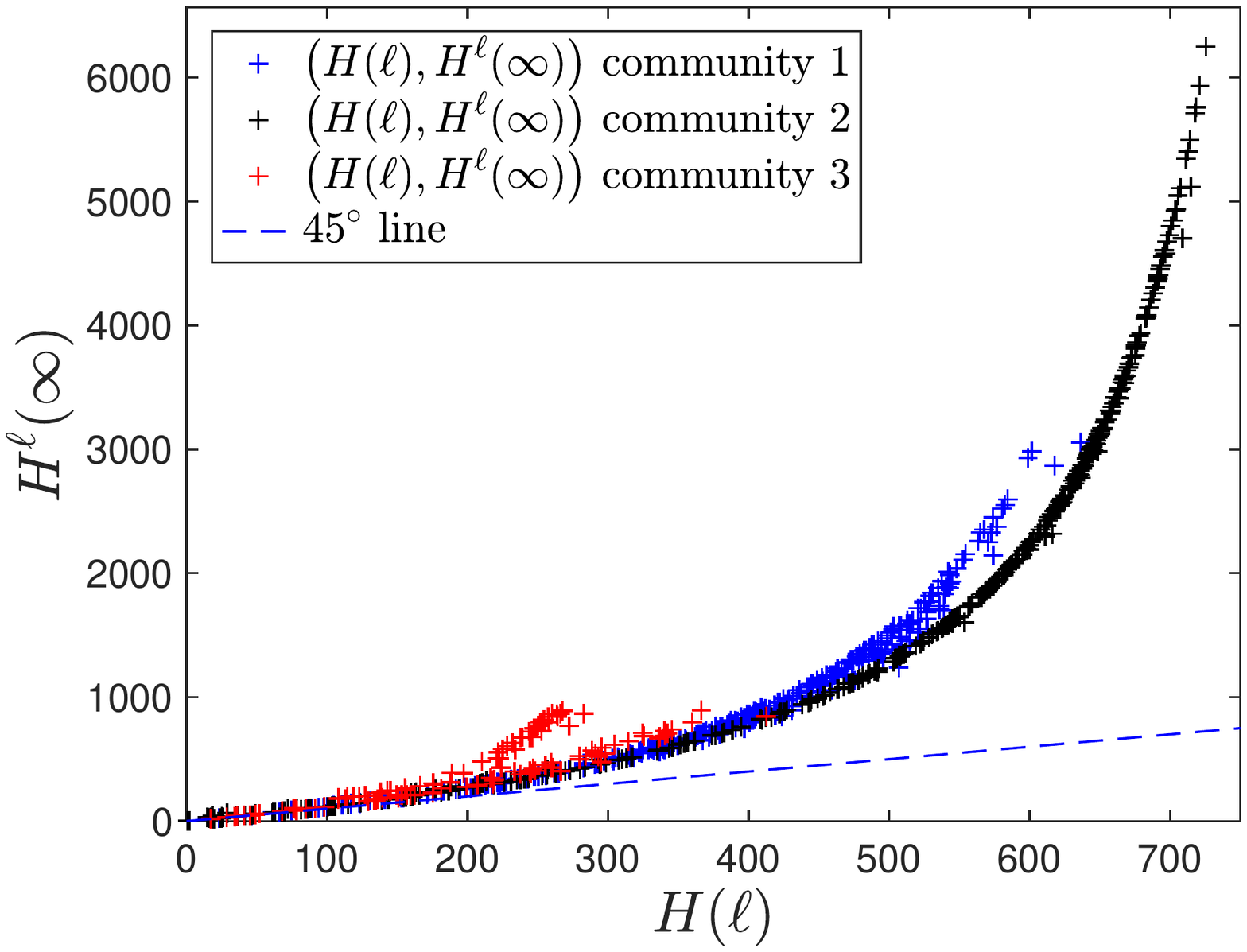}%
	}
	\caption{\label{fig:comm_mm_1_no_seed_HH_col} The asymptotic values $H^\ell(\infty)$ of the harmonic influence computed by the MPA against the corresponding exact values $H(\ell)$ computed by the definition, for the graph $\GG_1$. The different colors distinguish the crosses corresponding to nodes of the three different communities. All crosses are above the $45^\circ$ line.  }
\end{figure}

The simulation on the graph $\GG_1$ presents interesting differences. 
Fig.~\ref{fig:comm_mm_1_no_seed_dyn} represents the convergence of the algorithm. The $\Wij$ messages take about 15 steps to converge while $H^\ell(t)$ take about 34000 steps. 
Since graph $\GG_1$ is larger than $\GG_2$, but less dense, this moderate decrease in the convergence time is consistent with the $O(m/n)$ conjecture.
Fig.~\ref{fig:comm_mm_1_no_seed_HH_col} compares the estimates of the harmonic influence $\GG_1$ with their exact counterpart; the crosses are colored according to the community of the nodes. The black crosses correspond to nodes in the (largest) second community, the blue crosses to nodes in the first community and the red ones to nodes in the third (and smallest) community. 
The community structure produces an interesting artefact, which is made apparent by the alignment of the crosses: the MPA assigns excess influence to leaders within smaller communities, compared to leaders of larger communities.

\section{The Effects of Network Changes}\label{sec:T-VAR}

The structure of the MPA makes it easily adaptable to networks that change while the distributed computations are unfolding. New nodes and links may appear while other might disappear: if these changes happen without notice, the MPA cannot be restarted from the proper initial condition~\eqref{eq:MPA-initial}. For the sake of this discussion, we assume that the network changes only once, after the MPA has reached convergence on the initial network.
We may interpret the dynamics that follows the change as an MPA running on the new network, but starting from with a different initialization. The result in~\cite{RF:2016:TNSE} does not guarantee the convergence, because its proof uses a monotonicity property of the dynamics of $\Wij(t)$ that is only valid for the standard initialization. 
Moreover, it is not clear whether the MPA would carry some memory of the initial network.

In order to test these two facts, we construct a pair of small networks with partly similar topologies but very different harmonic influence profiles. The networks are  $\GG_3 = (I_3,E_3,C_3)$ and $\GG_4 = (I_4,E_4,C_4)$, with $I_3 = I_4 = \{\field, 1,\ldots, 50\}\,,$
while the edge sets $E_3$ and $E_4$ have the form of a ``wheel'' with additional connections. 
Both sets contain all the possible edges involving the field node, i.e.\ $\{\field, i\}$ for $i\in \{1,\ldots,50\},$
and the cycle
$$\{\{1,2\}, \{2,3\}, \ldots, \{49,50\}, \{1,50\}\} $$
connecting among all non-field nodes. 
Both sets also contain some extra edges of the form $\{i,j\}$ with $i\notin\{\field,1,26\}$ and $j\neq \field$: these are included with probability $p = 0.01$. 
Up to here, sets $E_3$ and $E_4$ are identical. To distinguish the networks, we include some additional edges $\{1, j\}$ in $E_3$ and some additional edges $\{26, j\}$  in $E_4$: we pick these edges with probability $q = 0.25$.
The matrix $C_3\in \Rp^{I_3\times I_3}$ has entries
\begin{equation*}
\left\{ \begin{array}{ll} 
	(C_3)_{i\field}= (C_3)_{\field i} = 0.040  & \text{ for every } i \in \{1,\ldots,50\} \\ 
	(C_3)_{ij}= 1 & \text{ if } i,j\neq \field \text{ and } \{i,j\}\in E_3 \\
	(C_3)_{ij}= 0 & \text{ if } \{i,j\}\notin E_1 \end{array}\right. 
\end{equation*}
The entries of $C_4$ are chosen similarly. 
The MPA starts on the network $\GG_3$ and, after a sufficiently large time $T$, 
continues on the network $\GG_4$. We use a bar to denote messages and estimates in this scenario. During the network change, the messages $\bWij(T)$ and $\bHij(T)$ corresponding to edges in $E_3\cap E_4$ retain their values, while the messages corresponding to missing edges in $E_3\setminus E_4$ are simply dropped. 
The messages corresponding to new edges $\{i,j\}\in E_4\setminus E_3$ are initialized by
$$\bWij(T) = 1\,,\bHij(T) = 1 \text{ where } i\neq \field\,. $$

\begin{figure}	
\centering 	
\fbox{%
	\includegraphics[trim={\figtrimlFB} {\figtrimbFBX} {\figtrimrFB} {\figtrimtFBX},
	clip, width={\figwidth}, keepaspectratio=true]{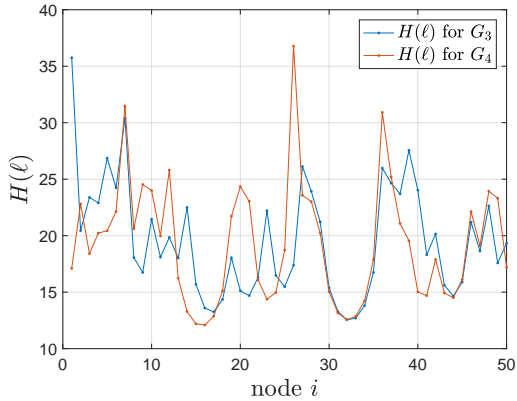}%
	}
	\caption{\label{fig:HH_esatti_G3G4} The harmonic influence $H(\ell)$ of the nodes of $\GG_3$ compared with the harmonic influence of the corresponding node in $\GG_4$.}
\end{figure}

We have repeated the simulation multiple times finding consistent results; we discuss one of the outcomes in what follows.
The exact profiles of the harmonic influence are compared in Fig.~\ref{fig:HH_esatti_G3G4}. The most influential nodes are node~1 in $\GG_3$ and node~26 in $\GG_4$; their influences change significantly between the two graphs. Some other nodes  hold  very similar influences in $\GG_3$ and $\GG_4$, e.g. nodes from 28 to 34. 

%
%

%

\begin{figure}	
\centering 	
\fbox{%
	\includegraphics[trim={\figtrimlFB} {\figtrimbFBX} {\figtrimrFB} {\figtrimtFBX},
	clip, width={\figwidth}, keepaspectratio=true]{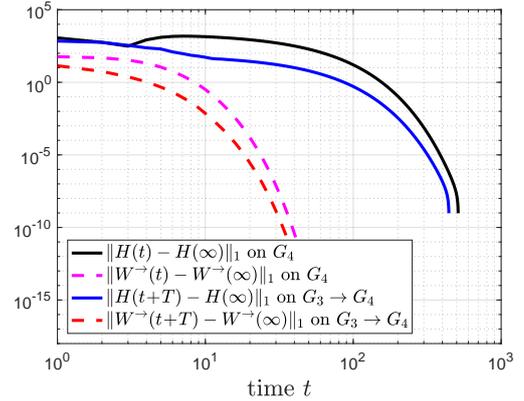}%
	}
	\caption{\label{fig:conv_time_G3G4} 
	The convergence of the MPA on the changing network scenario (after the change of topology) compared with the convergence over the network $\GG_4$. After the change of topology the MPA requires less iterations that a complete restart.	}
\end{figure}

We compare the convergence and estimates of the MPA on the changing scenario, with the convergence and estimate of the MPA started directly on $\GG_4$. The MPA on the changing network converges and requires less additional iterations (after $T$) to converge than the MPA on $\GG_4$: see Fig.~\ref{fig:conv_time_G3G4}. 
The convergence values $\bWij(\infty)$ and $\bar H^\ell(\infty)$ on the changing network coincide exactly with those obtained directly on $\GG_4$. 
This result led us to conclude that the convergence values only depend on the final topology.  

In conclusion, the MPA appears to be able to adapt smoothly to unforeseen changes in the network topology. Mathematically, we conjecture that it converges under general initial conditions and that it has 
a unique equilibrium. While giving a full proof of this conjecture remains an open problem, we have so far verified that the equilibrium $\Wij(\infty)$ is locally asymptotically stable. 

%



\end{document}